\documentclass[12pt]{article}

\usepackage{amsmath,amssymb}

\def \1{{\bf 1}}
\def \a{{{\mathfrak a}}}

\def \A{{\mathbb A}}

\def \bs{\backslash}

\def \C{{\mathbb C}}

\def \CO{{\cal O}}

\def \det{{\rm det}}

\def \eqn{\begin{eqnarray*}}
\def \endeqn{\end{eqnarray*}}
\def \eps{\epsilon}
\def \F{{\mathbb F}}

\def \Ga{\Gamma}

\def \GL{{\rm GL}}

\def \H{{\mathbb H}}

\def \Im{{\rm Im}}

\def \N{\mathbb N}

\def \p{{{\mathfrak p}}}

\def \qed{\ifmmode\eqno Q.E.D.\else\noproof\vskip 12pt plus 3pt minus 9pt \fi}
 \def\noproof{{\unskip\nobreak\hfill\penalty50\hskip2em\hbox{}%
     \nobreak\hfill Q.E.D.\parfillskip=0pt%
     \finalhyphendemerits=0\par}}
\def \Q{\mathbb Q}

\def \R{{\mathbb R}}
\def \Re{{\rm Re \hspace{1pt}}}

\def \ra{\rightarrow}

\def \SL{{\rm SL}}
\def \Spec{{\rm Spec}\ }

\def \Z{\mathbb Z}
\def \={\ =\ }

\newcommand{\rez}[1]{\frac{1}{#1}}

\newcounter{lemma}
\newcounter{corollary}
\newcounter{proposition}
\newcounter{theorem}

\newtheorem{theorex}{\hspace{-30pt}\stepcounter{conjecture} \stepcounter{lemma}
    \stepcounter{corollary} \stepcounter{proposition}Theorem}[section]

\begin{document}

\pagestyle{myheadings} \markright{PANORAMA OF ZETA
FUNCTIONS}

\title{Panorama of zeta functions}
\author{Anton Deitmar}

\date{}
\maketitle \tableofcontents

\begin{center}
{\bf Introduction}
\end{center}
In this essay I will give a strictly subjective selection
of different types of zeta functions. Instead of providing
a complete list, I will rather try to give the central
concepts and ideas underlying the theory.

Talking about zeta functions in general one inevitably is
led to start with the Riemann zeta function $\zeta(s)$. It
is defined as a {\it Dirichlet series:}
$$
\zeta(s)\=\sum_{n=1}^\infty n^{-s},
$$
which converges for each complex number $s$ of real part
greater than one. In the same region it possesses a
representation as a {\it Mellin integral:}
$$
\zeta(s)\=\rez{\Ga(s)}\int_0^\infty \rez{e^t-1}\ t^s
\frac{dt}t.
$$
This integral representation can be used to show that
$\zeta(s)$ extends to a meromorphic function on the
complex plane and satisfies the {\it functional equation:}
$$
\hat\zeta(s)\=\hat\zeta(1-s),
$$
where $\hat\zeta$ is the completed zeta function
$$
\hat\zeta(s)\= \pi^{-\frac s2}\Ga(\frac s2)\zeta(s).
$$
The zeta function has a simple pole at $s=1$ and is
regular otherwise. It has zeros at $-2,-4,\dots$ which are
called the trivial zeros. All other zeros lie in the strip
$0<\Re(s)<1$ and the famous {\it Riemann Hypothesis}
states that they should all lie at $\Re(s)=\rez 2$. This
conjecture was stated in the middle of the 19th century and has not
been proved to this day.\\
The Riemann Hypothesis is by no means the only riddle posed
by the Riemann zeta function. There is, for example, the
question about the spacing distribution of consecutive
zeros: Let $(\rho_n)_{n\in\N}$ be the ascending sequence
of imaginary parts of the zeros of $\zeta(s)$ in $\{
\Im(s)>0\}$, and let $\tilde\rho_n=\frac{\rho_n\log
\rho_n}{2\pi}$ be the normalized sequence. Let
$\delta_n=\tilde\rho_{n+1}-\tilde\rho_n$ be the sequence of
normalized spacings. Computations of pair correlation
functions \cite{Montgomery} and extensive numerical
calculations \cite{Odlyzko} then lead one to expect that
for any ``nice'' function $f$ on $(0,\infty)$
$$
\rez N\sum_{n=1}^N f(\delta_n)\ \ra\ \int_0^\infty f(s)
P(s) ds,\ \ \ {\rm as}\ N\ra\infty.
$$
Here $P(s)$ is the spacing distribution of a large random
Hermitian matrix. The function $P(s)$ vanishes to second
order at $s=0$. So, unlike a Poisson process, the numbers
$\tilde\rho_n$ ``repel'' each other.\\
This expectation is known as the ``GUE hypothesis'', where
$GUE$ stands for Gaussian Unitary Ensemble and describes
the spacing function $P$. The GUE hypothesis is widely
accepted, but far from proven.

Using the uniqueness of the prime decomposition of a
natural number the Dirichlet series can be turned into an
{\it Euler product:}
$$
\zeta(s)\=\prod_p(1-p^{-s})^{-1},
$$
which is extended over all prime numbers $p$ and converges
for real part of $s$ greater than one. The Euler product
indicates that the zeta function may be viewed as a means
to encode the infinite set of data given by the prime
numbers into a single object, the meromorphic function
$\zeta$. Indeed, the analytic behaviour of the zeta
function is exploited in the proof of the prime number
theorem \cite{ingham}, which says that the counting
function of prime numbers
$$
\pi(x)\=\#\{ p\le x,\ p\ {\rm prime}\}
$$
has the asymptotic behaviour
$$
\pi(x)\ \sim\ \frac x{\log x}
$$
as $x$ tends to $\infty$.

More generally the term zeta function is used for
generating series which encode infinite sets of data such
as the numbers of solutions  of an algebraic equation over
finite fields or the lengths of closed geodesics on a
Riemann surface. A zeta function is usually given as a
Dirichlet series, an Euler product or a Mellin integral.

Erich K\"ahler was very interested in zeta functions and
he even defined a new one himself. To explain this
consider the Euler product of the Riemann zeta function.
The set of prime numbers can be identified with the set of
all nonzero prime ideals of the ring of integers $\Z$, i.e.
the set of all closed points of the scheme $\Spec\Z$. This
is the starting point of a generalization to an arbitrary
arithmetic scheme $X$, whose so called {\it Hasse-Weil
zeta function} \cite{serre} is defined as a product over
the closed points of $X$:
$$
\zeta_X(s)\=\prod_{x\in |X|}(1-N(x)^{-s})^{-1},
$$
where for $x$ a closed point, $N(x)$ denotes the
order of the residue class field $\kappa(x)$ of $x$. \\
On the other hand, the Dirichlet series of $\zeta$ can be
interpreted as a sum over all nonzero ideals of $\Z$,
i.e., a sum over all nontrivial closed subsets of
$\Spec\Z$. In general, the Hasse-Weil zeta function of an
arithmetic scheme $X$ is also expressible as a Dirichlet
series, but not one that runs over all nontrivial closed
subsets of $X$. So K\"ahler defined a new zeta function in
the following way: For every closed subscheme $Y$ of $X$ he
defined $N(Y)=\prod_{y\in Y} N(y)$, and then
$Z_X(s)=\sum_{Y} N(Y)^{-s}$, where $N(Y)^{-s}$ is
considered to be zero if $N(Y)$ is infinite. In his thesis
\cite{Lu}, Lustig proves that this series converges for
real part of $s$ bigger than $2$, provided that the
dimension of $X$ does not exceed $2$. In his Monadologie
\cite{Kae-Mo}, 271-331, K\"ahler determined the zeta
function of an arithmetic curve of degree 2. His
calculation is reproduced in a modern language in
\cite{Berndt}. It is shown that the zeta function in this
case is a product of the Riemann zeta functions and a
Dirichlet $L$-series with an entire function whose zeros
can be given explicitly and lie at $\Re(s)=\rez 2$. If the
dimension of $X$ exceeds $2$, Witt (\cite{Witt} p.369)
showed that the series diverges for every complex number
$s$.

In the sequel, I will discuss three types of zeta, or
$L$-functions. The latter are a slight generalization of
zeta functions: they arise, for example by allowing twists
by characters in the coefficients of a Dirichlet series.
The first section will be on zeta and $L$-functions of
arithmetic origin which also might be called ``of
algebraic geometric origin'', since they are attached to
objects of algebraic geometry over the integers.\\
In the second section we will consider automorphic
$L$-functions. They are defined in an analytic setting and
look very different from the arithmetic ones of the
previous section, however, it is widely believed, and
proven in a number of cases, that automorphic
$L$-functions and
arithmetic $L$-functions are basically the same.\\
The third section will be devoted to zeta functions of
differential geometric origin. These are attached to
objects from differential geometry or global analysis,
like closed geodesics in Riemannian manifolds.

\section{Zeta and $L$-functions of arithmetic origin}
This section will be concerned with zeta functions whose
defining data come from number theory. Since the Riemann
zeta function has an Euler product over the primes, one
may view the prime numbers as defining data for $\zeta$,
so the Riemann zeta function falls into this category.

Let $K$ be a number field, i.e., a finite extension of the
field of rational numbers. The {\it Dedekind zeta
function} of $K$ is defined as
$$
\zeta_K(s)\=\sum_\a N(\a)^{-s},\ \ \ \Re(s)>1,
$$
where the sum runs over all nontrivial ideals $\a$ of the
ring of integers $\CO_K$ of $K$, and $N(\a)=\# (\CO_K/\a)$
is the norm of the ideal $\a$. The function $\zeta_K(s)$
can also be expressed as an Euler product
$$
\zeta_K(s)\=\prod_\p (1-N(\p)^{-s})^{-1}
$$
which runs over all nontrivial prime ideals $\p$ of
$\CO_K$. Further, $\zeta_K(s)$ extends to a meromorphic
function on the entire plane satisfying a functional
equation \cite{neukirch}. There are several
generalizations of these zeta functions, called
$L$-functions, which mostly are obtained by ``twisting by
characters''. For example
$$
L(\chi,s)\= \prod_{\p}(1-\chi(\p)N(\p)^{-s})^{-1}
$$
is the Hecke $L$-function attached to a character $\chi$
of the group $I_K$ of all fractional ideals of $K$ (see
\cite{neukirch}).

In his thesis \cite{tate}, J. Tate used harmonic analysis
on the adele ring $\A=\A_K$ and the idele group
$\A^\times=\A_K^\times$ to give quite far reaching
generalizations of the theory of $L$-functions which we
will now explain. A rational differential form $\omega$
will give a Haar measure on $\A^\times$, which, by the
product formula, being independent of the choice of
$\omega$, is called the {\it Tamagawa measure} on
$\A^\times$. The continuous homomorphisms of $\A^\times
/K^\times$ to $\C^\times$, also called {\it
quasi-characters}, are all of the form $z\mapsto
c(z)=\chi(z)|z|^s$ for a character $\chi$ and a complex
number $s$. For such a quasi-character and a sufficiently
regular function $f$ on $\A$ the general zeta function of
Tate is defined by
$$
\zeta(f,c)\=\int_{\A^\times}f(z) c(z) d^\times z,
$$
the integration being with respect to the Tamagawa measure
and convergent if $\Re(s)>1$. Let $\hat f$ be the additive
Fourier transform of $f$. Considered as an analytic
function in $c$ the zeta function is regular except for
two simple poles at $c_0$ and $c_1$, where $c_0(z)=1$ and
$c_1(z)=|z|$, with residues $-\kappa f(0)$ and $\kappa\hat
f(0)$ respectively, where
$$
\kappa\=\frac{2^{r_1}(2\pi)^{r_2}hR}{w\sqrt d}.
$$
Here $r_1$ and $r_2$ are the numbers of real embeddings,
resp. pairs of complex embeddings of $K$, $h$ is the class
number, $R$ the regulator, $w$ the order of the group of
roots of unity in $K$, and $d$ is the absolute value of the
discriminant. \\
Furthermore, there is a functional equation
$$
\zeta(f,c)\=\zeta(\hat f,\hat c),
$$
where $\hat c(z)=|z|c(z)^{-1}$.

One recovers classical zeta and $L$-functions by
restricting $c$ and $f$ to special cases. For example
$c(z)=|z|^s$ gives, for suitable $f$, the (completed)
Dedekind zeta function of $K$.

For $f$ of simple type the integral over $\A^\times$
defining $\zeta(f,c)$ can be written as a product of
integrals over the completions of $K$, so $\zeta(f,c)$
becomes a product of {\it local zeta functions} which
themselves satisfy local functional equations and bear
other interesting information.\\
J. Igusa \cite{igusa} replaced the characters in Tate's
local zeta functions by characters composed with
polynomials and so defined a new rich class of local zeta
functions with applications to other areas of number
theory.

Perhaps the most important generalization of the Riemann
zeta function is the class of motivic $L$-functions
\cite{deligne}. To explain this let us go back to the
Hasse-Weil zeta function
$\zeta_X(s)=\prod_{x\in|X|}(1-N(x))^{-1}$ of an arithmetic
scheme $X$. If $X$ is proper and flat over $\Spec\Z$ with
smooth generic fibre $X\otimes\Q$, then it will have good
reduction at almost all primes $p$, we call such primes
``good''. The Lefschetz trace formula for $l$-adic
cohomology implies that the zeta function equals a finite
number of Euler factors multiplied by
$$
\prod_{\nu=0}^{2\dim X} \prod_{p\ \rm good}
\det_{\Q_l}(1-p^{-s} Fr_p^* |
H^\nu(X\otimes\overline{\Q}_p,\Q_l))^{(-1)^{\nu+1}}.
$$
The inner product runs over all but finitely many primes.
Here $l=l_p$ is a prime different from $p$ and $Fr_p$ is
the geometric Frobenius at $p$. The characteristic
polynomial is known to have coefficients in $\Q$ which are
independent of $l$ (for $p$ good). It has turned out in a
number of cases that the individual factors
$$
\prod_{p\ \rm good} \det_{\Q_l}(1-p^{-s} Fr_p^* |
H^\nu(X\otimes\overline{\Q}_p,\Q_l))^{(-1)^{\nu+1}}
$$
themselves have a meromorphic continuation to $\C$, and
indeed satisfy a functional equation if suitably completed
at the bad primes. This then is considered to be the
$L$-function attached to the motive $H^\nu(X)$. Motives
form a conjectural category in which a scheme $X$
decomposes into $H^\nu(X)$ for $\nu=0,\dots,2\dim X$ and
all usual cohomology theories factor over this category.
This category is supposed to be large enough to contain
twists and all $L$-functions mentioned so far can be
realized as motivic $L$-functions.

There are various conjectures which relate vanishing orders
or special values of motivic $L$-functions to other
arithmetic quantities. We will here only give one example,
the conjecture of Birch and Swinnerton-Dyer. To explain
this let $E$ be an elliptic curve over a number field,
i.e., a projective curve of genus one with a fixed
rational point. Then there is a natural structure of an
abelian algebraic group on $E$. The group $E(K)$ of
rational points is known to be finitely generated, so its
{\it rank}, which is defined by $r=\dim_\Q E(K)\otimes\Q$,
is finite. The Birch and Swinnerton-Dyer conjecture states
that $r$ should be equal to the vanishing order of the
Hasse-Weil zeta function $\zeta_E(s)$ at $s=1$. There is
also a more refined version \cite{Tate} of this conjecture
giving an arithmetical interpretation of the first
nontrivial Taylor coefficient of $\zeta_E(s)$ at $s=1$.

\section{Automorphic $L$-functions}
Let $f$ be a cusp form of weight $2k$ for some natural
number $k$ as in \cite{serre}, i.e., the function $f$ is
holomorphic on the upper half plane $\H$ in $\C$, and has a
certain invariance property under the action of the modular
group $\SL_2(\Z)$ on $\H$. Then $f$ admits a Fourier
expansion
$$
f(z)\=\sum_{n=1}^\infty a_n e^{2\pi i zn}.
$$
Define its $L$-function for $\Re(s)>1$ by
$$
L(f,s)\= \sum_{n=1}^\infty \frac{a_n}{n^s}.
$$
The easily established integral representation
$$
\hat L(f,s)\= (2\pi)^{-s}\Ga(s)L(f,s)\=\int_0^\infty f(it)
t^s \frac{dt}t,
$$
implies that $L(f,s)$ extends to an entire function
satisfying the functional equation $\hat L(f,s)=(-1)^k\hat
L(f,2k-s)$. With $\Lambda(f,s)=\hat L(f,2ks)$ this becomes
$$
\Lambda(f,s)\= (-1)^k \Lambda(f,1-s).
$$
This construction can be extended to cusp forms for
suitable subgroups of the modular group. These
$L$-functions look like purely analytical objects with no
connection to the $L$-functions of arithmetic origin
mentioned earlier. Thus it was particularly daring of A.
Weil, G. Shimura, and Y. Taniyama  in 1955 to propose the
conjecture that the zeta function of any elliptic curve
over $\Q$ coincides with a $\Lambda(f,s)$ for a suitable
cusp form $f$. This conjecture was proved in part by A.
Wiles and R. Taylor \cite{Wiles, Wiles-Taylor}  providing a
proof of Fermat's Last Theorem as a consequence.
Subsequently, the conjecture has been proved in full by
Breuil, Conrad, Diamond and Taylor \cite{Breuiletal}.

The upper half plane is a homogeneous space of the group
$\SL_2(\R)$, and so cusp forms may be viewed as functions
on this group, in particular, they are vectors in the
natural unitary representation of $\SL_2(\R)$ on the space
$$
L^2(\SL_2(\Z)\bs \SL_2(\R)).
$$
Going even further one can extend this quotient space to
the quotient of the adele group $\GL_2(\A)$ modulo its
discrete subgroup $\GL_2(\Q)$, so cusp forms become
vectors in
$$
L^2(\GL_2(\Q)\bs \GL_2(\A)^1),
$$
where $\GL_2(\A)^1$ denotes the set of all matrices in
$\GL_2(\A)$ whose determinant has absolute value one. Now
$\GL_2$ can be replaced by $\GL_n$ for $n\in\N$ and one can
imitate the methods of Tate's thesis (the case $n=1$) to
arrive at a much more general definition of an automorphic
$L$-function: this is an Euler product $L(\pi,s)$ attached
to an automorphic representation $\pi$ of $\GL_n(\A)^1$,
i.e. an irreducible subrepresentation $\pi$ of
$L^2(\GL_n(\Q)\bs \GL_n(\A)^1)$. As in the $\GL_1$-case it
has an integral representation as a Mellin transform and it
extends to a meromorphic representation, which is entire
if $\pi$ is cuspidal and $n>1$. Furthermore it satisfies a
functional equation
$$
L(\pi,s)\=\eps(\pi,s) L(\tilde\pi,1-s),
$$
where $\tilde\pi$ is the contragredient representation and
$\eps(\pi,s)$ is a constant multiplied by an exponential
\cite{GodementJacqet}.

Extending the Weil-Shimura-Taniyama conjecture, R.P.
Langlands conjectured  in the 1960s that any motivic
$L$-function coincides with $L(\pi,s)$ for some cuspidal
$\pi$. An affirmative answer to this question would prove
many other, older conjectures in number theory. It has
local and characteristic $p$ analogues, which have been
proved \cite{Henniart,Lafforgue}.

\section{Zeta functions of differential-geometric\\ origin}

In this section we will be concerned with zeta functions
attached to objects arising in differential geometry. These
will serve as a measure of complexity of the geometric
objects. We start with an oversimplified example: a finite
graph $X$. If it has no closed paths, it is topologically
trivial, so the set of primitive closed paths gives a
measure of complexity. Primitivity here means that a path
(or a walk, as some people say) $c_0$ is not the power of a
shorter one, so you walk each closed path only once. The
number of closed paths will in general be infinite, so one
considers the following zeta function as a formal power
series at first:
$$
Z_X(T)\= \prod_{c_0}(1-T^{l(c_0)}),
$$
where the product runs over all primitive closed paths and
$l(c_0)$ denotes the length of a given path $c_0$. Then it
turns out \cite{Bass} that $Z_X(T)$ is in fact a
polynomial, for it can be written as $\det(1-TA)$ for some
generalized adjacency operator $A$ for the graph $X$.

If one replaces the graph $X$ by a compact Riemannian
surface $Y$ of genus $\ge 2$, then one can attach a natural
hyperbolic metric to $Y$ and replace the paths by closed
geodesics. One ends up with the {\it Ruelle Zeta Function}
$$
R_Y(s)\=\prod_c(1-e^{-sl(c)}),
$$
where the product runs over all primitive closed geodesics
in $Y$ and $l(c)>0$ is the length of the geodesic $c$. The
Ruelle Zeta Function $R_Y(s)$ equals $Z_Y(s)/Z_Y(s+1)$,
where $Z_Y(s)$ is the {\it Selberg Zeta Function} attached
to $Y$, defined by
$$
Z_Y(s)\=\prod_c\prod_{n\ge 0} (1-e^{(s+n)l(c)}).
$$
The Selberg Zeta Function can be studied using harmonic
analysis and one can prove that it extends to an entire
function having all its zeros in the set $\R\cup (\rez
2+i\R)$, i.e., $Z_Y(s)$ satisfies a generalized Riemann
Hypothesis \cite{Hejhal}.

Note that in this section, Euler products often occur
without the $(-1)$ in the exponent. This is more than a
matter of taste, since, for example, the Selberg zeta
function is entire this way round. The reason why one has
to take different exponents becomes transparent when one
generalizes the Selberg zeta function to spaces of higher
rank \cite{geom}: the natural exponents are certain Euler
characteristics which can take positive or negative values.

The set of closed geodesics on $Y$ is in bijection with
the set of closed orbits of the geodesic flow $\phi$ on
the sphere bundle $SY$. So $R_Y(s)=R_\phi(s)$ is a special
case of a {\it dynamical zeta function}, since it counts
closed orbits of a dynamical system. Actually D. Ruelle
proved the meromorphicity of $R_\phi(s)$ in a more general
setting: he needed $\phi$ only to be a smooth flow
satisfying a certain hyperbolicity condition. He used the
theory of Markov families to express $R_\phi(s)$ as a
quotient of certain transfer operators. For the geodesic
flow of the modular curve $\SL_2(\Z)\bs \H$ the
eigenvectors of these transfer operators can be identified
with modular functions \cite{Lewis}.

Also for discrete dynamical systems there is a theory of
zeta functions. Let most generally $f$ be an invertible
self map of a set $X$ and define its zeta function as the
formal power series
$$
Z_f(T)\=\prod_{p} (1-T^{l(p)}),
$$
where the product runs over the set of periodic orbits and
$l(p)\in\N$ is the period of $p$. One has to put
restrictions on $f$ in order for $Z_f$ to be well defined.
For $f$ being a diffeomorphism of a compact manifold
satisfying certain natural hyperbolicity conditions A.
Manning has shown in \cite{Manning}, that $Z_f(T)$ indeed
is a rational function. This in particular means that all
fixed point data of $f$ can be exhibited from the finite
set
of poles and zeros of $Z_f(T)$.\\
The Hasse-Weil zeta function of a smooth projective
variety $V$ over a finite field $\F_q$ with $q$ elements
can be viewed as a dynamical zeta function too. It
coincides with $Z_{Fr}(q^{-s})$, where $Fr$ is the
Frobenius acting on $V(\overline{\F_q})$.

One feature that makes zeta functions of geometric origin
so attractive is that they tend to satisfy {\it Lefschetz
formulae}. For example, as Artin, Grothendieck and Verdier
have shown, the Hasse-Weil zeta function of a smooth
projective variety $V$ over $\F_q$ satisfies
$$
\zeta_V(s)\=\prod_{\nu=0}^{2\dim V}\det(1-q^{-s}Fr |
H^\nu(V,\Q_l))^{(-1)^{\nu+1}}.
$$
The Selberg zeta function satisfies a similar Lefschetz
formula involving the Frobenius vector field and the
tangential cohomology of the contracting foliation
\cite{Guillemin,Patterson,geom}. This fact has given rise
to some far reaching conjectures whether these formulae are
valid for more general systems and if they even could be
part of a cohomological framework which would explain most
of the conjectural properties of zeta and $L$-functions of
arithmetic schemes \cite{Deninger}.

\section{Closing remarks}
In this brief survey we have missed out many other
important classes of zeta functions. It was not our aim to
give an exhaustive list but to show some general lines of
development. Summarizing we find that zeta functions,
given by a Dirichlet series, an Euler product or a Mellin
integral, encode infinite sets of data. The first thing one
usually asks for, is convergence, next meromorphicity and
functional equation. If one is lucky the function turns
out to be rational, this then means that the infinite set
of data can be recovered from the finitely many poles and
zeros. One finally starts to ask in which way the analytic
behaviour of the zeta function reflects properties of the
encoded objects. Prominent examples of this are the Prime
Number Theorem, in which the position of the pole and the
zeros of the zeta function give information on the growth
of the data, or results or conjectures on special values
like the Birch and Swinnerton-Dyer conjecture.

Seeking harmony and simplicity the human mind is always
tempted to believe that objects of similar behaviour,
although of very different origin, should be of the same
nature. Remarkably enough this has turned out true in the
case of the Taniyama-Shimura-Weil conjecture and has shown
strong evidence in the case of the Langlands conjecture.

Whenever some entities are counted with some mathematical
structure on them, it is likely that a zeta function can
be set up and often enough it will extend to a meromorphic
function. Zeta functions show up in all areas of
mathematics and they encode properties of the entities
counted which are well hidden and hard to come by
otherwise. They easily give fuel for bold new conjectures
and thus drive on the progress of mathematics. It is a
fairly safe assertion to say that zeta functions of
various kinds will stay in the focus of mathematical
attention for times to come.

{\scriptsize University of Exeter, Mathematics, Exeter EX4
4QE, England}


\small
\begin{thebibliography}{XXX}

\bibitem{Bass}
 \bf Bass, H.:
 \it The Ihara-Selberg zeta function of a tree lattice.
 \rm Int. J. Math. 3, No.6, 717-797 (1992).

\bibitem{Berndt}
 \bf Berndt, R.:
 \it K\"ahler's computation of his Zeta
 function for an arithmetic curve of degree two.
 \rm Mitt. Math. Ges. Hamburg III. Hamburg 1985.

\bibitem{Breuiletal}
\bf Breuil, C.; Conrad, B.; Diamond, F.; Taylor, R.:
 \it On the modularity of elliptic curves over $\Q$: wild
3-adic exercises.
 \rm J. Amer. Math. Soc. 14,  843-939 (2001).

\bibitem{geom}
 \bf Deitmar, A.:
 \it Geometric zeta-functions of locally symmetric spaces.
 \rm Am. J. Math. 122, vol 5, 887-926 (2000).

\bibitem{deligne}
 \bf Deligne, P.:
 \it Valeurs de fonctions $L$ et périodes d'intégrales.
 \rm  With an appendix by N. Koblitz and A. Ogus. Proc.
Sympos. Pure Math., XXXIII, Automorphic forms,
representations and $L$-functions (Proc. Sympos. Pure
Math., Oregon State Univ., Corvallis, Ore., 1977), Part 2,
pp. 313--346, Amer. Math. Soc., Providence, R.I., 1979.

\bibitem{Deninger}
 \bf Deninger, C.:
 \it Some analogies between number theory and dynamical systems on foliated spaces.
 \rm Proceedings of the International Congress of
Mathematicians, Vol. I (Berlin, 1998). Doc. Math., Extra
Vol. I, 163-186 (1998).

\bibitem{GodementJacqet}
 \bf Godement, R.; Jacquet, H.:
 \it Zeta functions of simple algebras.
 \rm Lecture Notes in Mathematics. 260. Berlin-Heidelberg-New York: Springer-Verlag. (1972).

\bibitem{Guillemin}
 \bf Guillemin, V.:
 \it Lectures on spectral theory of elliptic operators.
 \rm Duke Math. J. 44, 485-517 (1977).

\bibitem{Hejhal}
 \bf Hejhal, D.:
 \it The Selberg trace formula for $PSL_2(\R)$ I.
 \rm Springer Lecture Notes 548, 1976.

\bibitem{Henniart}
 \bf Henniart, G.:
 \it Une preuve simple des conjectures de Langlands pour
 ${\rm GL}(n)$ sur un corps $p$-adique.
 \rm Invent. Math. 139, 439-455 (2000).

\bibitem{igusa}
 \bf Igusa, J.:
 \it An introduction to the theory of local zeta functions.
 \rm AMS/IP Studies in Advanced Mathematics,
14. American Mathematical Society, Providence, RI;
International Press, Cambridge, MA, 2000.

\bibitem{ingham}
 \bf Ingham, A. E.:
 \it The distribution of prime numbers.
 \rm Reprint of the 1932 original. With a foreword by R. C.
Vaughan. Cambridge Mathematical Library. Cambridge
University Press, Cambridge, 1990.

\bibitem{Kae-Mo}
 \bf K\"ahler, E.:
 \it Monadologie III.
 \rm Hamburg 1985.

\bibitem{Lafforgue}
 \bf Lafforgue, L.:
 \it Chtoucas de Drinfeld et correspondance de Langlands.
 \rm Invent. math.; DOI 10.1007/s002220100174 (2001).

\bibitem{Lewis}
 \bf Lewis, J.; Zagier, D.:
 \it Period functions and the Selberg zeta function for the modular group.
 \rm The mathematical beauty of physics (Saclay, 1996), 83--97,
Adv. Ser. Math. Phys., 24, World Sci. Publishing, River
Edge, NJ, 1997.

\bibitem{Lu}
 \bf Lustig, G:
 \it \"Uber die Zetafunktion einer arithmetischen
 Mannigfaltigkeit.
 \rm Math. Nachr. 14, 309-330 (1955).

\bibitem{Manning}
 \bf Manning, A.:
 \it Axiom A diffeomorphisms have rational zeta functions.
 \rm Bull. Lond. Math. Soc. 3, 215-220 (1971).

\bibitem{Montgomery}
 \bf Montgomery, H. L.:
 \it The pair correlation of zeros of the zeta function.
 \rm Analytic number theory (Proc. Sympos. Pure
Math., Vol. XXIV, St. Louis Univ., St. Louis, Mo., 1972),
pp. 181-193. Amer. Math. Soc., Providence, R.I., (1973).


\bibitem{neukirch}
\bf Neukirch, J.: \it Algebraic number theory. \rm
Grundlehren der Mathematischen Wissenschaften. 322.
Berlin: Springer 1999.

\bibitem{Odlyzko}
 \bf Odlyzko, A. M.:
 \it On the distribution of spacings between zeros of the zeta function.
 \rm Math. Comp. 48, no.177, 273-308 (1987).


\bibitem{Patterson}
 \bf  Patterson, S. J.:
 \it On Ruelle's zeta-function.
 \rm Festschrift in honor of I. I. Piatetski-Shapiro on the occasion
of his sixtieth birthday, Part II (Ramat Aviv, 1989),
163--184, Israel Math. Conf. Proc., 3, Weizmann,
Jerusalem, 1990.

\bibitem{Ruelle}
 \bf Ruelle, D.:
 \it Zeta functions for Expanding maps and Anosov flows.
 \rm Invent. math. 34, 231-244 (1976).

\bibitem{sel}
 \bf Selberg, A:
 \it Harmonic Analysis and Discontinuous Groups in weakly symmetric Riemannian spaces with Applications to Dirichlet Series.
 \rm J. Indian. Math. Soc. 20, 47-87 (1956).

\bibitem{serre}
 \bf Serre, J.-P.:
 \it Zeta and L functions.
 \rm Arithmetical algebraic Geom., Proc. Conf. Purdue Univ. 1963, 82-92 (1965).

\bibitem{tate}
 \bf Tate, J.:
 \it Fourier analysis in number fields, and Hecke's zeta-functions.
 \rm Algebraic Number Theory (Proc.
Instructional Conf., Brighton, 1965) Thompson, Washington,
D.C. 305-347 (1967).

\bibitem{Tate}
 \bf Tate, J.:
 \it On the conjectures of Birch and
Swinnerton-Dyer and a geometric analog.
 \rm Séminaire
Bourbaki, Vol. 9, Exp. No. 306, 415-440, Soc. Math.
France, Paris, 1995.

\bibitem{Wiles}
 \bf Wiles, A:
 \it Modular elliptic curves and Fermat's last theorem.
 \rm Ann. of Math. (2) 141, 443-551 (1995).

\bibitem{Wiles-Taylor}
 \bf Taylor, R.; Wiles, A.:
 \it Ring-theoretic properties of certain Hecke algebras.
 \rm Ann. of Math. (2) 141, 553-572 (1995).

\bibitem{Witt}
 \bf Witt, E.:
 \it Collected papers. Gesammelte Abhandlungen. (German)
\rm With an essay by Guenter Harder on
Witt vectors. Edited and with a preface in English and
German by Ina Kersten. Springer-Verlag, Berlin, 1998.

\end{thebibliography}
\end{document}